\documentclass[a4paper,12pt]{article}
   \usepackage[top=2.5cm,bottom=2.5cm,left=2.5cm,right=2.5cm]{geometry}
    \usepackage{amssymb}
    \usepackage{graphicx}
    \usepackage{amsthm}
    \usepackage{amsmath}
    \usepackage[english]{babel}
    \usepackage[latin1]{inputenc}
    \usepackage{textcomp}
    \usepackage{longtable}
    \usepackage{listings}

\newtheorem{theorem}{Theorem}

\newtheorem{lemma}{Lemma}

\begin{document}

\title{\vspace*{2.5cm} New Bounds for the Sum of Powers of Normalized Laplacian Eigenvalues of Graphs}

\author{{Gian Paolo Clemente$^{a}$}, Alessandra Cornaro$^{a}$}
\date{}
\maketitle

\vspace*{-5mm}

\begin{center}
$^{{}^a}$\emph{Department of Mathematics and Econometrics, Catholic
University, Milan, Italy}\\[2mm]
{\tt
gianpaolo.clemente@unicatt.it, alessandra.cornaro@unicatt.it}\\[5mm]
\vspace{5mm}

\end{center}

\vspace{3mm}

\thispagestyle{empty}
\begin{abstract}
For a simple and connected graph, a new graph invariant $s_{\alpha}^{*}(G)$, defined as the sum of $\alpha$ powers of the eigenvalues of the normalized Laplacian matrix, has been introduced by Bozkurt and Bozkurt in \cite{Bozkurt2012}. Lower and upper bounds have been proposed by the authors.
In this paper, we localize the eigenvalues of the normalized Laplacian matrix by adapting a theoretical method,
proposed in Bianchi and Torriero (\cite{BT}), based on majorization techniques.
Through this approach we derive upper and lower bounds of $s_{\alpha}^{*}(G)$.
Some numerical examples show how sharper results can be obtained with respect to those existing in literature.
\end{abstract}
\baselineskip=0.30in

\section{Introduction}\label{intro}

Among the various indices in Mathematical Chemistry, a whole new family of descriptors $s_{\alpha}^{*}(G)$, defined as the sum of $\alpha$ powers of the eigenvalues of the normalized Laplacian matrix, has been proposed by Bozkurt and Bozkurt in \cite{Bozkurt2012}.
These authors found a number of bounds for arbitrary $\alpha$ and particularly for  $\alpha=-1$, which is the case of the degree Kirchhoff Index.
Recently, Bianchi et al. proposed a variety of lower and upper bounds for $s_{\alpha}^{*}(G)$ in \cite{BCPT2} and for the Kirchhoff Index in \cite{BCPT1} derived via majorization techniques. In particular, the authors showed that it is possible to obtain tighter results taking into account additional information on the localization of the eigenvalues of proper matrices associated to the graph. From a theoretical point of view, some well--known inequalities on the localization of real eigenvalues have been provided in literature and they can be used to compute the above mentioned bounds.\\
Alternative inequalities involving the localization of some eigenvalues of the transition matrix of the graph have been numerically computed in \cite{CC1} and \cite{CC} by adapting a theoretical methodology proposed in Bianchi and Torriero \cite{BT} based on nonlinear global optimization problems solved through majorization techniques. By means of these results, tighter lower bounds for the Kirchhoff Index for some classes of graphs have been derived in \cite{CC1}.\\
The purpose of this paper is to use this fruitful theoretical method providing indeed some formulae that allow us to compute lower bounds for the first and the second eigenvalues of the normalized Laplacian matrix in a fairly straightforward way.
Thus, we obtain new bounds for $s_{\alpha}^{*}(G)$ considering both non-bipartite and bipartite graphs.

In Section 2 some preliminaries are given. In Section 3 we describe the nonlinear optimization problem useful for our analysis. Solving this optimization problem, lower bounds of the first and second eigenvalue of the normalized Laplacian matrix have been provided in Section 4. Finally, in Section 5 several numerical results are reported showing how the proposed bounds are tighter than those given in \cite{Bozkurt2012}.

\section{Notations and preliminaries}\label{not}

In this section we first recall some basic notions on graph theory. For more details refer to \cite{Wilson}.

\noindent Let $G=(V,E)$ be a simple, connected, undirected graph where $V=\{1, 2, \ldots, n\}$ is the set of vertices and $E\subseteq V\times V$ the set of edges, $|E|=m$.

\noindent The degree sequence of $G$ is denoted by $\pi=(d_{1},d_{2},..,d_{n})$ and it is arranged in non-increasing order $d_{1}\geq d_{2}\geq \cdots \geq d_{n}$, where $d_{i}$ is the degree of
vertex $i$.
It is well known that $\overset{n}{\underset{%
i=1}{\sum }}d_{i}=2m$ and that if $G$ is a tree, i.e. a connected graph
without cycles, $m=n-1.$ Let $A(G)$ be the adjacency matrix of $G$ and $D(G)$
be the diagonal matrix of vertex degrees. The matrix $L(G)=D(G)-A(G)$ is
called Laplacian matrix of $G$, while $\mathcal{L}%
(G)=D(G)^{-1/2}L(G)D(G)^{-1/2}$ is known as normalized Laplacian matrix. Let $\mu _{1}\geq \mu _{2}\geq ...\geq \mu
_{n}$
 be the (real) eigenvalues of $L(G)$ and $\lambda _{1}\geq \lambda _{2}\geq ...\geq \lambda _{n}$ be the (real) eigenvalues of $\mathcal{L}(G)$. The  following properties of spectra of $L(G)$ and $\mathcal{L}(G)$ hold:
\begin{equation*}
\overset{n}{\underset{i=1}{\sum }}\mu _{i}=\mathrm{tr}(L(G))=2m;\,\,\,\mu
_{1}\geq 1+d_{1}\geq \dfrac{2m}{n};\,\,\,\mu _{n}=0,\text{ }\mu _{n-1}>0.
\end{equation*}%
\begin{equation*}
\overset{n}{\underset{i=1}{\sum }}\lambda _{i}=\mathrm{tr}(\mathcal{L}%
(G))=n;\,\,\,\overset{n}{\underset{i=1}{\sum }}\lambda _{i}^{2}=\mathrm{tr}(%
\mathcal{L}^{2}(G))=n+2\sum_{(i,j)\in E}\frac{1}{d_{i}d_{j}};\,\,\lambda
_{n}=0;\lambda _{1}\leq 2.
\end{equation*}%

\noindent Our aim is the analysis of a particular topological index, $s_{\alpha}^{*}(G)$.
Topological indices have been widely explored in different fields, i.e. mathematical chemistry and more recently in complex network analysis.
In particular, Zhou (see \cite{Zhou}) proposed the index:
\begin{equation*}
s_{\alpha }(G)=\sum_{i=1}^{n-1}\mu _{i}^{\alpha },\alpha \neq 0,1,
\end{equation*}
defined as the sum of the $\alpha$-th power of the non-zero Laplacian eigenvalues of a graph $G$.

Over the last years this index and its bounds have been intensely studied:
Zhou (see \cite{Zhou}) established some properties of $s_\alpha(G)$ and some improvements have been provided in \cite{LiuLiu}, \cite{Tian},  \cite{Zhou2} and \cite{ZhouIlic}.
In \cite{BCT2}, taking into account the Schur-convexity or Schur-concavity of the functions $s_{\alpha}(G)$ for $\alpha >1$ and $\alpha <0$ or $0<\alpha <1$ respectively, the same bounds as in \cite{Zhou} have been derived. Furthermore, considering additional information on the localization of the eigenvalues,
the authors provide also sharper bounds.

Bozkurt and Bozkurt in \cite{Bozkurt2012} introduced parallely to \cite{Zhou} the following new graph invariant:
\begin{equation*}
s_{\alpha }^{\ast }(G)=\sum_{i=1}^{n-1}\lambda _{i}^{\alpha },\alpha \neq 0,1,
\end{equation*}
characterized as the sum of the $\alpha$-th power of the non-zero normalized Laplacian eigenvalues of a graph.
Several properties of this index have been proposed in \cite{Bozkurt2012}
and some lower and upper bounds for a connected graph have been derived.

In \cite{BCPT2}, considering the Schur-convexity or Schur-concavity of the functions $s_{\alpha }^{\ast }(G)$ and using additional information on the localization of the eigenvalues, the following Theorems, which generalize Theorem 3.3 in \cite{Bozkurt2012}, have been proved.

\begin{theorem}
\label{th:lambda_1}
Let $G$ be a simple connected graph with $n\geq 3$\ vertices and $\lambda _{1} \geq \theta$:
\begin{enumerate}
\item if $\alpha <0$ or $\alpha >1$ then
\begin{equation}\label{eq:simple1}
s^*_{\alpha}(G) \ge \theta^{\alpha}+ \frac {(n-\theta)^{\alpha}} {(n-2)^{\alpha-1}}
\end{equation}

\item if $0<\alpha <1$ then
\begin{equation}\label{eq:simple2}
s_{\alpha }^{\ast }(G)\leq \theta^{\alpha }+\frac{(n-\theta)^{\alpha }}{(n-2)^{\alpha
-1}}.
\end{equation}
\end{enumerate}
\end{theorem}

\begin{theorem}
\label{th:lambda_2}
Let $G$ be a simple connected graph with $n\geq 4$\ vertices which is not complete and $\lambda _{1} \geq \theta$, $\lambda_{2}\geq \beta $ with $\theta\geq\beta$ and $\theta+\beta(n-2) > n$.
\begin{enumerate}
\item if $\alpha <0$ or $\alpha >1$ then
\begin{equation}
s_{\alpha }^{\ast }(G)\geq
\theta^{\alpha }+\beta  ^{\alpha }+\frac{(n-\theta-\beta)^{\alpha }}{%
(n-3)^{\alpha -1}} \label{eq:simpleb1}
\end{equation}

\item if $0<\alpha <1$ then
\begin{equation}s_{\alpha }^{\ast }(G)\leq \theta^{\alpha }+
\beta^{\alpha }+\frac{(n-\theta-\beta)^{\alpha }}{(n-3)^{\alpha -1}}. \label{eq:simpleb2}
\end{equation}
\end{enumerate}
\end{theorem}

It is noteworthy to state that the results in Theorem \ref{th:lambda_2} are tighter than those in Theorem \ref{th:lambda_1} (for more details see \cite{BCT2} and \cite{BCT1}).

In \cite{Bozkurt2012}, the bounds in Theorem \ref{th:lambda_1} has been previously proved identifying $\theta$ as
\begin{equation*}
P=1+\sqrt{\frac{2}{n(n-1)}\sum_{(i,j \in E)}\frac{1}{d_{i}d_{j}}}.
\end{equation*}

In Section \ref{sec:nb}, by adapting some results proved in \cite{BT}, we provide lower bounds of $\lambda_{1}$ and $\lambda_{2}$ that enable us to obtain tighter bounds than in \cite{Bozkurt2012}. In what follows we refer to the values of lower bounds of $\lambda_{1}$ and $\lambda_{2}$ as $Q$ and $R$, respectively, i.e. $\lambda_{1}\geq Q$ and $\lambda_{2}\geq R$.

\section{A Nonlinear optimization problem to bound eigenvalues}\label{sec:majo}
We now recall a methodology based on majorization techniques (see \cite{BT} and \cite{Marshall}) that allows us to find a suitable localization
of $\lambda_{1}$ and $\lambda_{2}$ in order to provide some new results for bounds of $s^*_{\alpha}(G)$.

\noindent At this regard, we define the set
\[S^{\lambda}_{b}=\{\boldsymbol{\lambda}\in \mathbb{R}_{+}^{n-1}: \lambda_{1} \geq \lambda_{2}\geq ... \geq \lambda_{n-1}\geq 0,
\sum_{i=1}^{n-1}\lambda_{i}=n,g(\boldsymbol{\lambda})=\sum_{i=1}^{n-1}\lambda_{i}^{p}=b\},\]
where $p$ is an integer greater than 1.


\noindent  The following fundamental lemma holds (see Lemma 2.1 in \cite{BT}):
\begin{lemma}
Fix $b \in \mathbb{R}$ and consider the set $S^{\lambda}_{b}$.
Then either $b=\frac{n^{p}}{(n-1)^{p-1}}$ or there exists a unique integer $1 \leq h^{*} < (n-1)$ such that:
\begin{equation}
\frac{n^{p}}{\left( h^{\ast }+1\right) ^{p-1}}<b\leq \frac{n^{p}}{\left(h^{\ast }\right) ^{p-1}},
\label{eq:cond}
\end{equation}
where $  h^{\ast}=\left\lfloor \sqrt[p-1]{\frac{n^{p}}{b}}\right\rfloor .$ \\
\label{le:lem}
\end{lemma}

The following Theorem is derived by Theorem 3.2 in \cite{BT}. Moreover, a lower bound for $\lambda_{h}$ by solving the following optimization problem:
\begin{center}
\ \ \ \ \ \ \ \ \ \ \ \ \ \ \ \ \ \ \ \ \ \ \ \ min $(\lambda_{h})$ \ subject to $\boldsymbol{\lambda} \in S^{\lambda}_{b}$ \ \ \ \ \ \ \ \ \ \ \ \ \ \ \ \ \ \ \ \ \ \ \ \ \ \ \ \ $P^{*}(h)$
\end{center}

\begin{theorem}
The solution of the optimization problem $P^{*}(h)$ is
$(\frac{n}{n-1})$ if $b=\frac{n^{p}}{(n-1)^{p-1}}$. \\

\noindent If $b \neq \frac{n^{p}}{(n-1)^{p-1}}$, the solution of the optimization problem $P^{*}(h)$ is $\delta^{*}$ where
\noindent \begin{enumerate}
\item {for $h=1$, $\delta^{*}$ is the unique root  of the equation
\begin{equation}
f(\delta,p) =  h^{\ast }\delta^{p}+ (n-h^{\ast }\delta)^{p}- b= 0
\label{eq:f3}
\end{equation}
in $I=\left( \frac{n}{h^{\ast }+1},\frac{n}{h^{\ast }}\right]$;
}

\item {for $1 < h \leq (h^{*}+1)$, $\delta^{*}$ is the unique root  of the equation

\begin{equation}
f(\delta,p) = (n-h)\delta^{p}+ (h-1)\frac{(n-(n-h)\delta)^{p} }{(h-1)^{p}} - b=0
\label{eq:fn}
\end{equation}
in $I=(0,\frac{n}{n-1}]$;
}
\item {for $h >(h^{*}+1)$, $\delta^{*}$  is zero. }

\end{enumerate}
\label{th:teo1}
\end{theorem}

\section{New bounds for normalized Laplacian eigenvalues} \label{sec:nb}

In order to find new bounds for $s^{\ast}_{\alpha}(G)$, we make use of the methodology introduced in Section \ref{sec:majo}
that allows us to provide lower bounds for $\lambda_{1}$ and $\lambda_{2}$.
At this regard, we now consider Theorem \ref{th:teo1} limiting the analysis when $p=2$:
in this case we know indeed that $b=n+2 \sum_{(i,j) \in E} \frac{1}{d_{i}d{j}}$.

When $b=\frac{n^{2}}{(n-1)}$, then the solution of optimization problem is $\frac{n}{(n-1)}$.
This is the case of the complete graph $K_{n}$.
Instead, when $b \neq \frac{n^{2}}{(n-1)}$, $h^{\ast}=\left\lfloor {\frac{n^{2}}{b}}\right\rfloor.$

Since we want to get a lower bound for $\lambda_{1}$, we solve equation (\ref{eq:f3}).
The acceptable solution in the interval $I$ is equal to ${\delta}^{\ast}=\frac{\left(n+\sqrt{\frac{b(h^{\ast}+1)-n^{2}}{h^{\ast}}}\right)}{(1+h^{\ast})}$ and we refer to this value as $Q$.

The value of $Q$ can be compared to $P$ in order to show how bounds (\ref{eq:simple1}) and (\ref{eq:simple2}), computed by assuming $\theta=Q$, perform better than those for $\theta=P$.

It is well known that, for every connected graph of order $n$ (see \cite{BCPT1}), we have:
\begin{equation}
\frac{1}{n-1} \leq \frac{2}{n} \sum_{(i,j)\in E}\frac{1}{d_{i}d_{j}} <1,
\label{eq:h}
\end{equation}
and the left inequality is attained for the complete graph $G=K_{n}$.

%
%

\begin{figure}[!h]
\centering
\includegraphics[scale=0.4]{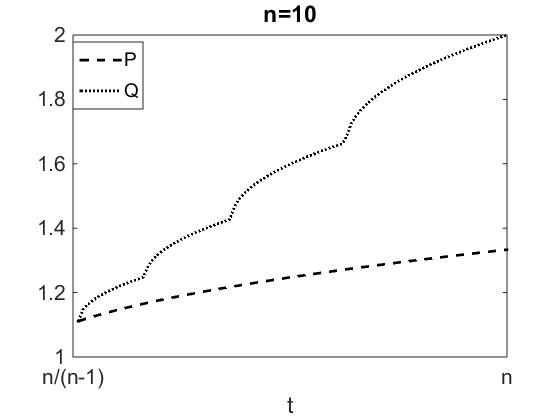}
\includegraphics[scale=0.4]{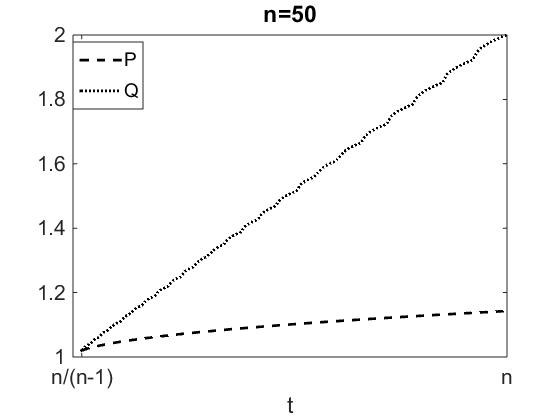}
\includegraphics[scale=0.4]{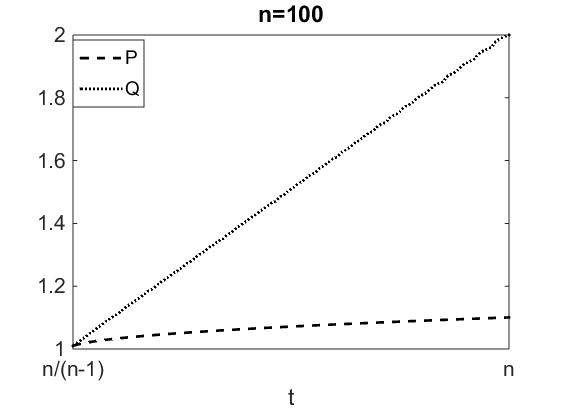}
\includegraphics[scale=0.4]{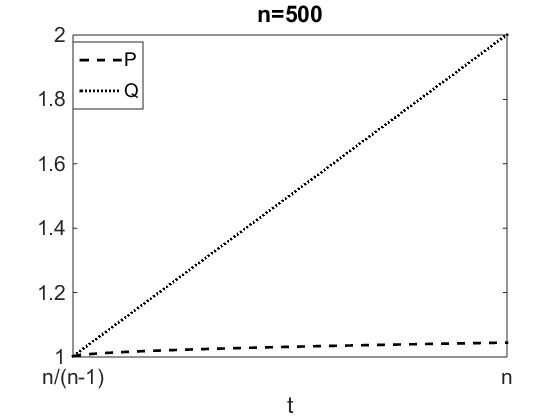}
\caption{$Q$ and $P$ according to different values of $t \in \left(\frac{n}{n-1},n\right)$ and several number of vertices.}
\label{fig:FigureQP1}
\end{figure}

Figure \ref{fig:FigureQP1} reports patterns of $P$ and $Q$, varying the quantity $t=2 \cdot \sum_{(i,j)\in E}\frac{1}{d_{i}d_{j}}$ in the proper interval $\left(\frac{n}{n-1},n\right)$ (see Equation (\ref{eq:h})) for alternative values of number of vertices $n$.

Being $P=1+\sqrt{\frac{t}{n(n-1)}}$, it is easy to see that $P$ has a monotonic behaviour with respect to $t$ in the interval $\left(\frac{n}{n-1},1+\sqrt{\frac{1}{n-1}}\right)$.
Figure \ref{fig:FigureQP1} shows that $Q$ assumes values equal to $P$ when $t = \frac{n}{n-1}$ (indeed $P=Q=\frac{n}{n-1}$ in this case) and then it increases faster than $P$ in the considered interval.

\begin{figure}[!h]
\centering
\includegraphics[scale=0.2]{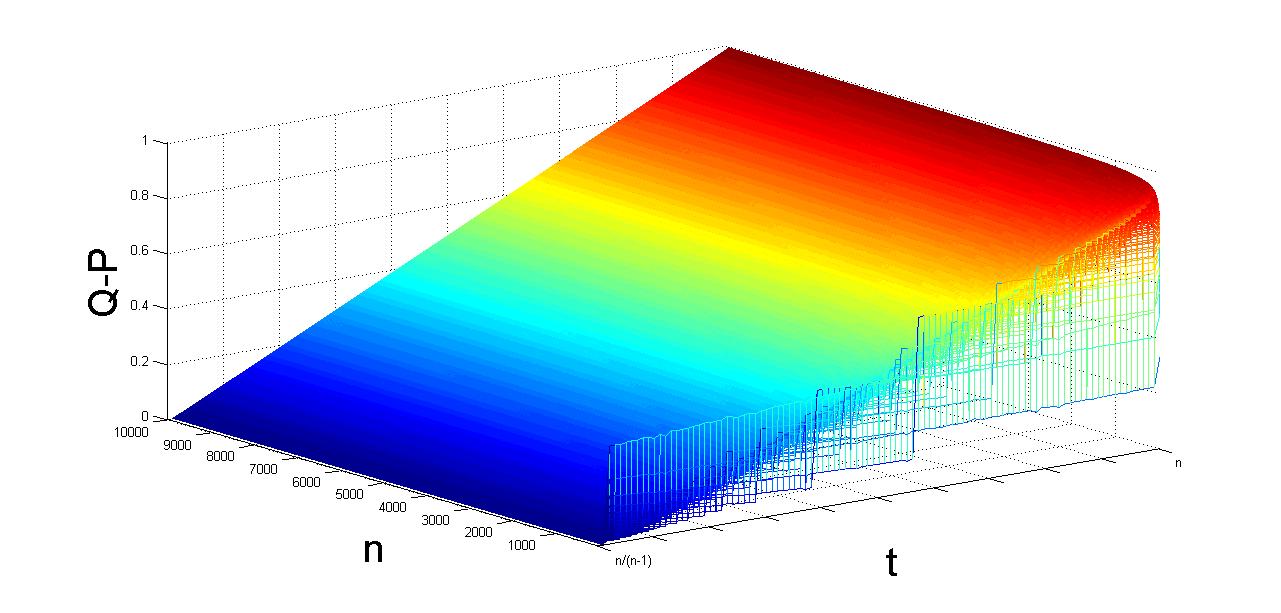}
\caption{$Q-P$ according to different $n$ and for values of $t \in \left(\frac{n}{n-1},n\right)$.}
\label{fig:FigureQP}
\end{figure}

Going deeply into the analysis, we graphically show that the difference $Q-P$ is always strictly greater than zero according to several values of $n$ and $t$ (see Figure \ref{fig:FigureQP}).

Since $\lambda_{1}\geq Q\geq P$, we provide bounds (\ref{eq:simple1}) and (\ref{eq:simple2}) better than in \cite{Bozkurt2012} (see \cite{BCT2} and \cite{BCT1} for more theoretical details).

With the aim to improve previous results, we can now derive additional information on $\lambda_{2}$.
To evaluate the lower bound $R$ of $\lambda_{2}$,
from Theorem \ref{th:teo1}, considering the case $h=2$,
since $h\leq (h^{*}+1)$, we solve the equation (\ref{eq:fn})
finding the acceptable solution:
\begin{equation*}
R=\delta^{\ast}=\frac{n-\sqrt{\frac{b(n-1)-n^{2}}{n-2}}}{n-1}.
\end{equation*}

In this case, the rightmost inequality in (\ref{eq:h}) implies $t < n$ and then $b < 2n$. By plugging this information in the value of $R$, we easily obtain $R \leq \frac{n}{n-1}$ that fulfills the condition $R \leq Q$ of Theorem \ref{th:lambda_2}. The other condition of Theorem \ref{th:lambda_2} , $Q+R(n-2) > n$, will be numerically checked in the next Section. In what follows, fixed values of $Q$ and $R$, we can numerically compute bounds (\ref{eq:simpleb1}) and (\ref{eq:simpleb2}) and comparing with those in \cite{Bozkurt2012}.

In case of bipartite graphs it is well known that $\lambda_{1}=2$.
If we set $\theta=2$ in Theorem \ref{th:lambda_1}, we derive the same results found in \cite{BCPT2}.
Furthermore, by placing $\theta=2$ and $\beta=R$ in bounds (\ref{eq:simpleb1}) and (\ref{eq:simpleb2}), we provide bounds also for bipartite graphs.

\section{Some Numerical Results}\label{sec:NumRes}

The proposed bounds have been evaluated on different graphs.
We now focus only on non-bipartite graphs and we provide a comparison with literature (see \cite{Bozkurt2012}).\\
In order to assure a robust analysis, graphs have been randomly generated following the Erd\"os-R\'enyi (ER) model $G_{ER}(n,q)$ (see \cite{Boll}, \cite{Chung}, \cite{ER59} and \cite{ER60}). Graphs have been obtained by using a MatLab code that gives back only connected graph based on the ER model (see \cite{CC1} and \cite{CC}). In this fashion, the graph is constructed by connecting nodes randomly such that edges are included with probability independent from every other edge. The results are based on a classic assumption of a probability of existence of edges $q$ equal to $0.5$. We obtain indeed that the generated graphs have a number of edges not far from the half of its maximum value as proved in the literature (see for example \cite{Estr}).

At this regard, in Table \ref{tab:res}, $s^{*}_{\alpha}(G)$ has been computed for several graphs by fixing $\alpha$ equal to 0.5.
We report values of upper bound (\ref{eq:simple2}) evaluated by using $\theta=Q$ or $\theta=P$ (as proposed in \cite{Bozkurt2012}).
We refer to these bounds as (\ref{eq:simple2}$Q$) and (\ref{eq:simple2}$P$).
Likewise bound (\ref{eq:simpleb2}$QR$) identifies bound (\ref{eq:simpleb2}) evaluated when  $\theta=Q$ and $\beta=R$, where the results has been provided assuring that assumptions of Theorem \ref{th:lambda_2} are satisfied.
Relative errors $r$ measures the absolute value of the difference between the upper bounds and $s^{*}_{\alpha}(G)$ divided by the value of $s^{*}_{\alpha}(G)$. We observe an improvement with respect to existing bounds according to all the analyzed graphs and the improvement appears reduced for very large graphs. However, for large graphs the formula provided in \cite{Bozkurt2012} already gives a very low relative error.

\begin{table}[!h]
\tiny
\centering
\begin{tabular}{|c|c|c||c|c|c|c||c|c|c||}
\hline\hline
$n$ & $d_{1}$ & $m$ & $s^{*}_{\alpha}(G)$ & bound (\ref{eq:simple2}$Q$) & bound (\ref{eq:simpleb2}$QR$) & bound (\ref{eq:simple2}$P$) & r(\ref{eq:simple2}$Q$) & r(\ref{eq:simpleb2}$QR$)  & r(\ref{eq:simple2}$P$) \\ \hline\hline

4	&2	&3	&3.35	&3.44	&3.43	&3.46	&2.86\%	&2.55\%	&3.47\% \\ \hline
5	&4	&9	&4.46	&4.47	&4.47	&4.47	&0.23\%	&0.21\%	&0.25\% \\ \hline
6	&3	&6	&5.30	&5.47	&5.46	&5.48	&3.13\%	&3.00\%	&3.27\% \\ \hline
7	&5	&14	&6.43	&6.48	&6.48	&6.48	&0.83\%	&0.81\%	&0.86\% \\ \hline
8	&5	&13	&7.33	&7.48	&7.48	&7.48	&2.02\%	&1.98\%	&2.06\% \\ \hline
9	&6	&16	&8.31	&8.48	&8.48	&8.48	&2.04\%	&2.01\%	&2.07\% \\ \hline
10	&8	&25	&9.39	&9.51	&9.48	&9.52	&1.36\%	&1.04\%	&1.37\% \\ \hline
20	&15	&95	&19.37	&19.51	&19.49	&19.51	&0.71\%	&0.62\%	&0.72\% \\ \hline
30	&19	&209	&29.36	&29.50	&29.50	&29.50	&0.49\%	&0.46\%	&0.49\% \\ \hline
50	&33	&604	&49.37	&49.50	&49.50	&49.50	&0.27\% &0.26\%	&0.27\% \\ \hline
100	&60	&2459	&99.37	&99.50	&99.50	&99.50	&0.13\%	&0.12\% 	&0.13\% \\ \hline
200	&116	&10001	&199.38	&199.50	&199.50	&199.50	&0.06\%	&0.05\%	&0.06\% \\ \hline
300	&179	&22437	&299.37	&299.50	&299.50	&299.50	&0.04\%	&0.04\%	&0.04\% \\ \hline
500	&279	&62456	&499.38	&499.50	&499.50	&499.50	&0.03\%	&0.02\%	&0.03\% \\ \hline \hline

\end{tabular}
\caption[]{Upper bounds for $s^{*}_{\alpha}(G)$ for $\alpha=0.5$ and relative errors.}
\label{tab:res}
\end{table}

\begin{figure}[!h]
\centering
\includegraphics[scale=0.6]{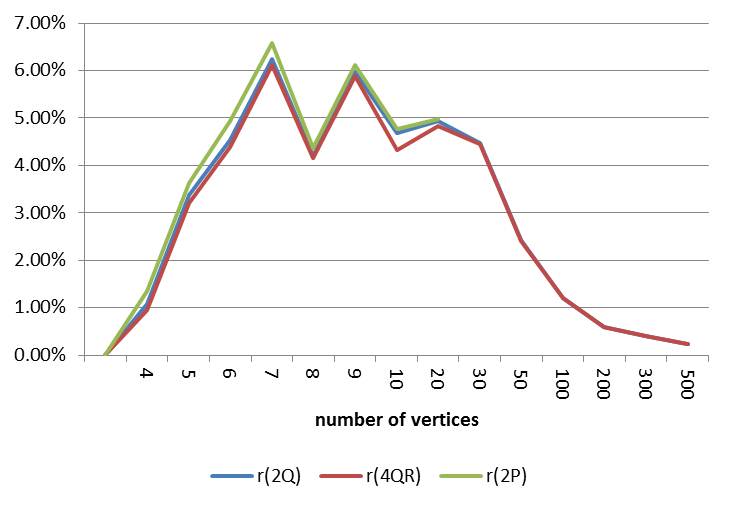}
\includegraphics[scale=0.6]{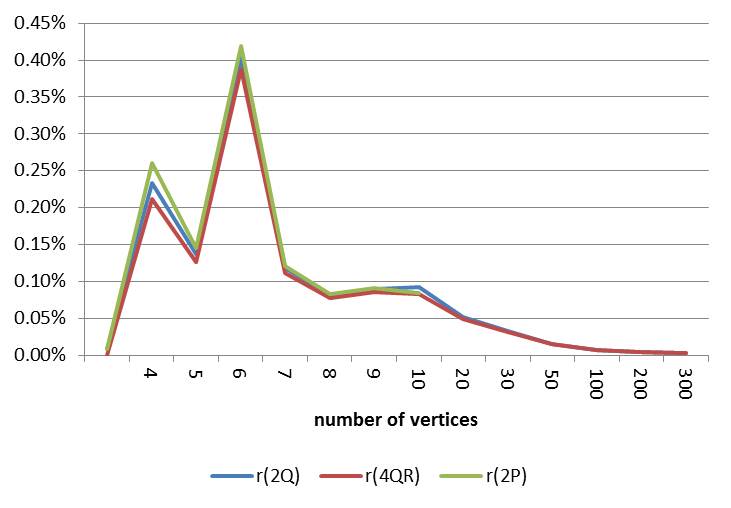}
\caption{Relative errors of upper bounds of $s^{*}_{0.5}(G)$ for graphs $ER(n,0.1)$ and $ER(n,0.9)$ respectively.}
\label{fig:F1}
\end{figure}

The comparison has been extended in order to test the behaviour of the upper bounds on alternative graphs.
First of all, in the ER model used to generate graphs, the parameter $q$ can be thought of as a weighting function. As $q$ increases from $0$ to $1$, the model becomes more and more likely to include graphs with more edges and less and less likely to include graphs with fewer edges.
In this regard, we assign several values of $q$ moving from the default value of $0.5$.
For sake of simplicity we report only the relative errors derived for graphs generated by using respectively $q=0.1$ and $q=0.9$ (see Figure \ref{fig:F1}).
In all cases bound (\ref{eq:simpleb2}$QR$) assures the best approximation to $s^{*}_{\alpha}(G)$ for $\alpha=0.5$.
We observe a best behaviour of all bounds when $q=0.9$ because we are moving towards the complete graph. We have indeed that the density of the graphs increases as long as greater probabilities are considered.

Finally, for the same index $s^{*}_{0.5}(G)$, upper bounds have been evaluated for trees\footnote{Tree has been generated by using Pr\"ufer code. The Pr\"ufer sequence of a labeled tree is a unique sequence associated to the tree. The sequence for a tree on $n$ vertices has length $n-2$ and it can be generated by a simple iterative algorithm. It is a way to map bijectively trees on $n$ vertices into $n-2$ long sequences of integers drawn from $n$.}.
Table \ref{tab:res2} depicts slighter differences for larger graphs in this case too. However it could be noticed how the relative improvement of bounds respect to other bounds is greater than in case of non-bipartite graphs.
Despite greater relative errors are observed, bound (\ref{eq:simpleb2}$QR$) is confirmed as the tighter bound also in this case.

\begin{table}[!h]
\tiny
\centering
\begin{tabular}{|c||c|c|c|c|}
\hline\hline
&\multicolumn{4}{c}{\textbf{Trees}} \vline  \\ \hline
$n$ & $s^{*}_{\alpha}(G)$ & r(\ref{eq:simple2}$Q$) & r(\ref{eq:simpleb2}$QR$)  & r(\ref{eq:simple2}$P$) \\ \hline\hline

4	 &3.35 	&2.04\%	&1.95\%	&3.47\%	\\ \hline
5	&4.32 	&2.17\%	&2.14\%	&3.48\%	\\ \hline
6	 &5.23 	&3.56\%	&3.51\%	&4.73\%	 \\ \hline
7	 &6.19 	&3.62\%	&3.59\%	&4.67\%	\\ \hline
8	 &7.15 	&3.67\%	&3.64\%	&4.61\%	\\ \hline
9	 &8.22 	&2.35\%	&2.34\%	&3.20\%	\\ \hline
10	 &8.85 	&6.35\%	&6.32\% &7.15\%	\\ \hline
20	 &18.07 	&7.45\%	&7.44\%	&7.88\%	\\ \hline
30	 &27.63 	&6.47\%	&6.47\%	&6.77\%	\\ \hline
50	 &45.73 	&8.07\%	&8.06\%	&8.25\%	\\ \hline
100	 &91.23 	&8.97\%	&8.97\%	&9.06\%	\\ \hline
200	 &182.72 	&9.14\%	&9.14\%	&9.18\%	\\ \hline
300	 &274.71 	&8.99\%	&8.99\%	&9.02\%	\\ \hline
500	 &457.71 	&9.11\%	&9.11\%	&9.13\% \\ \hline \hline

\end{tabular}
\caption[]{$s^{*}_{0.5}(G)$ and relative errors for Trees $T$.}
\label{tab:res2}
\end{table}

The analysis has been further developed considering a value of $\alpha$ equal to $1.5$.
Generating a similar sample of graphs, both $s^{*}_{1.5}(G)$ and the relative bounds have been derived.
For sake of simplicity we report only the results for $ER(n,0.5)$ observing that the additional information on the localization of $\lambda_{1}$ and $\lambda_{2}$ lead to the tighter lower bound (\ref{eq:simpleb1}$QR$).
Analogous results have been obtained by considering both ER graphs with alternative values of $q$ and bipartite graphs.

\begin{table}[!h]
\tiny
\centering
\begin{tabular}{|c|c|c||c|c|c|c||c|c|c||}
\hline\hline
$n$ & $d_{1}$ & $m$ & $s^{*}_{\alpha}(G)$ & bound (\ref{eq:simple1}$Q$) & bound (\ref{eq:simpleb1}$QR$) & bound (\ref{eq:simple1}$P$) & r(\ref{eq:simple1}$Q$) & r(\ref{eq:simpleb1}$QR$)  & r(\ref{eq:simple1}$P$) \\ \hline\hline
4	&3	&4	 &4.79 	 &4.66 	 &4.67 	 &4.62 	&2.69\%	&2.39\%	&3.49\% \\ \hline
5	&2	&4	 &6.22 	 &5.65 	 &5.69 	 &5.60 	&9.07\%	&8.51\%	&9.95\% \\ \hline
6	&4	&9	 &6.85 	 &6.59 	 &6.60 	 &6.57 	&3.78\%	&3.63\%	&4.00\% \\ \hline
7	&6	&13	 &7.77 	 &7.57 	 &7.57 	 &7.56 	&2.56\%	&2.49\%	&2.65\% \\ \hline
8	&7	&18	 &8.75 	 &8.56 	 &8.56 	 &8.55 	&2.15\%	&2.10\%	&2.20\% \\ \hline
9	&4	&12	 &10.28 	 &9.58 	 &9.59 	 &9.55 	&6.88\%	&6.79\%	&7.15\% \\ \hline
10	&7	&26	 &10.83 	 &10.52 	 &10.55 	 &10.51 	&2.90\%	&2.60\%	&2.94\% \\ \hline
20	&13	&98	 &20.88 	 &20.50 	 &20.52 	 &20.50 	&1.81\%	&1.71\%	&1.81\% \\ \hline
30	&19	&222	 &30.88 	 &30.50 	 &30.51 	 &30.50 	&1.21\%	&1.18\%	&1.21\% \\ \hline
50	&31	&644	 &50.85 	 &50.50 	 &50.51 	 &50.50 	&0.68\%	&0.67\%	&0.68\% \\ \hline
100	&62	&2512	 &100.87 	 &100.50 	 &100.50 	 &100.50 	&0.37\%	&0.36\%	&0.37\% \\ \hline
200	&117	&9918	 &200.88 	 &200.50 	 &200.50 	 &200.50 	&0.19\%	&0.19\%	&0.19\% \\ \hline
300	&179	&22540	 &300.87 	 &300.50 	 &300.50 	 &300.50 	&0.12\%	&0.12\%	&0.12\% \\ \hline
500	&279	&62063	 &500.88 	 &500.50 	 &500.50 	 &500.50 	&0.08\%	&0.08\%	&0.08\% \\ \hline \hline

\end{tabular}
\caption[]{Lower bounds for $s^{*}_{1.5}(G)$ and absolute value of relative errors.}
\label{tab:res3}
\end{table}

\newpage
\section{Conclusions}\label{Conc}
In this paper we provide tighter bounds for the sum of the $\alpha$-power of the non-zero normalized Laplacian eigenvalues taking into account additional information on the localization of the eigenvalues of the normalized Laplacian matrix of the graph, $\mathcal{L}(G)$.
To this aim lower bounds of the eigenvalues are derived by means of the solution of a class of suitable nonlinear optimization problems based on majorization techniques.
We provide indeed closed formulae that allow to compute upper and lower bounds of $s^{*}_{\alpha}(G)$ by using the additional information on the first and the second eigenvalue of $\mathcal{L}(G)$.
Numerical comparisons confirm how bounds are former than those existing in literature.
In particular, the analysis has been developed randomly generating both bipartite and non-bipartite graphs with a different number of vertices.

\bibliographystyle{plain}
\bibliography{biblio}

\end{document}